\documentclass[11pt]{amsart}
\usepackage{amsmath,amssymb,amsthm,upref,graphicx,mathrsfs}

\usepackage{color,xcolor}
\usepackage[
  colorlinks=true,
  linkcolor=blue,
  citecolor=blue,
  urlcolor=blue]{hyperref}

\numberwithin{equation}{section}

\newtheorem{thm}{Theorem}[section]
\newtheorem{cor}[thm]{Corollary}
\newtheorem{lem}[thm]{Lemma}

\newtheorem{df}{Definition}[section]

\numberwithin{equation}{section}

\begin{document}

%--------------------------------------------------
%%Don not change any thing in this part
\leftline{ \scriptsize}

\vspace{1.3 cm}
%----------------------------------------------------------------------------
\title
{Divisibility of LCM matrices by Totally nonnegative GCD matrices}
\author{Peeraphat Gatephan and Kijti Rodtes}
\thanks{{\scriptsize
}}
\hskip -0.4 true cm

\maketitle

%-----------------------------------------------------------------------------

%----------------------------------------------
\begin{abstract}  In this paper, we show that all totally nonnegative GCD matrices are always divisors of the corresponding LCM matrices in the ring $\mathbb{M}_{n}(\mathbb{Z})$.  We also introduce \lq\lq column monotone matrices" used to construct all totally nonegative GCD matrices. 
\end{abstract}

\vskip 0.2 true cm

%-----------------------------------------------------------------------------

\pagestyle{myheadings}
\markboth{\rightline {\scriptsize Peeraphat Gatephan and Kijti Rodtes}}
         {\leftline{\scriptsize }}
\bigskip
\bigskip

%-----------------------------------------------------------------------------
%-----------------------------------------------------------------------------

\vskip 0.4 true cm  
%--------------------------------------------\section{Introduction}------------------%--------------------------------------------\section{Introduction}------------------%--------------------------------------------\section{Introduction}------------------
\section{Introduction}
Let $S=\{x_1,\ldots,x_n\}$ be a set of distinct positive integers and $f$ be an arithmetical function. The $n\times n$ matrix $(S)=(s_{ij})$, where $s_{ij}=(x_{i},x_{j})$, the greatest common divisor of $x_{i}$ and $x_{j}$, is called the greatest common divisor matrix (GCD) matrix and $(f(S))$ denote the GCD matrix having $f$ evaluated at the $ij-$entry of GCD matrix $(S)$.  The least common multiple (LCM) matrix $[S]$ is defined similarly. The study of LCM matrices and GCD matrices began with the research of the famous number theorist H.J.S. Smith in 1876. Smith \cite{S1876} proposed that if $S=\{x_1,x_2,\ldots,x_n\}$ is factor-closed set, then the determinant of $\text{GCD matrix }(S)$, $\operatorname{det}(S)$, is $\prod\limits_{k=1}^n\phi(x_k)$, where $\phi(x_k)$ is the Euler's totient function.  A set $S$ is factor closed if it contains all divisors of $x$ for any $x\in S$. In the early years of the study, the researcher studied about the determinant of LCM matrices and GCD matrices on gcd-closed sets beginning with  Beslin and  Ligh in 1989. They \cite{BL1989} proposed the determinant of GCD matrices defined on gcd-closed set $S=\{x_1,x_2,\ldots,x_n\}$ which is \begin{eqnarray*}
	\operatorname{det}(S)=\prod\limits_{\substack{d\mid x_i,\,\,d\nmid x_t,\,\,t<i}}\phi(d).
\end{eqnarray*} 
A set $S$ is gcd-closed if $(x,y)\in S$ for all $x,y\in S$. Furthermore, they proposed that for any set of distinct positive integers $S$, $(S)$ is a positive definite matrix and they calculated the inverse formulas of LCM matrices and GCD matrices when their inverse exist. 	\\

For $A,B\in\mathbb{M}_{n}(\mathbb{Z})$ (the set of all $n\times n$ matrices over $\mathbb{Z}$), we say that $A$ divides $B$ or $B$ is divisible by $A$ in the ring $\mathbb{M}_{n}(\mathbb{Z})$ if there exists a matrix $C\in\mathbb{M}_{n}(\mathbb{Z})$ such that $B=AC$ or $B=CA$. We simply write $A\mid B$ if $A$ divides $B$ in ring $\mathbb{M}_{n}(\mathbb{Z})$ and $A\nmid B$ otherwise. Divisibility is an interesting topic in the study of LCM matrices and GCD matrices, starting with the research of Bourque and Ligh \cite{04} in 1992. They showed that if $S$ is a factor-closed set then $(S)\mid[S]$. Later in 2002, Hong \cite{13} showed that for any gcd-closed set $S$ with $\lvert S\lvert\leq3$, $(S)\mid[S]$ and there is a gcd-closed set $S$ with $\lvert S\rvert\geq4$ such that $(S)\nmid[S]$. In  the same paper, Hong raised the following open problem.\\\\
\textbf{Problem.\,}\cite{13}
Let $n\geq4$. Find a necessary and sufficient conditions on a gcd-closed set with $\lvert S\rvert=n$ such that $(S)\mid[S]$.\\

After that, many researchers studied the division of LCM matrices by GCD matrices on gcd-closed set. Nowadays, this problem has been answered in some special gcd-closed sets with $\lvert S\rvert=4$ by Zhao \cite{37} and $\lvert S\rvert=5$ by Zhao-Zhao \cite{38}. Further, the greatest-type divisor is also used to prove the problem, which was originally used by Hong. For $x,y\in S$ and $x<y$, if $x\mid y$ and the condition \lq\lq$x\mid z\mid y$ and $z\in S$ imply that $z\in\{x,y\}$" holds, then $x$ is said to be a greatest-type divisor of $y$ in $S$ and we denote by $G_{S}(y)$ the set of all greatest-type divisors of $y$ in $S$. In 2006, Hong \cite{16} showed that if $S$ is a gcd-closed set with $\operatorname{max}_{x\in S}\{\lvert G_{S}(x)\rvert\}=1$, then $(S)\mid[S]$. In addition to the aforementioned results, similar studies have been carried out by Feng, Hong and Zhao \cite{07} in cases where $S$ is a gcd-closed set satisfying $\operatorname{max}_{x\in S}\{\lvert G_{S}(x)\rvert\}\leq2$ and by  Altini\c{s}ik, Yildiz and Keskin \cite{00} in case where $S$ is a gcd-closed set satisfying $\operatorname{max}_{x\in S}\{\lvert G_{S}(x)\rvert\}\leq3$. On the other hand, some researchers paid attention to study divisibility of LCM matrices by GCD matrices in different sets, such as finite coprime divisor chain sets (see the definition in the discussion after Corollary \ref{co01}), see e.g. \cite{30,TLC2016,35}.\\

In the recent years, Guillot and Wu \cite{GWTN} proposed the necessary and sufficient conditions of totally nonnegative GCD matrices and a formula for the inverse of totally nonnegative GCD matrices. The matrix is totally nonnegative (TN) matrix if its all minors are nonnegative.  Note that some of them can not be constructed from any gcd-closed set nor set of finite many coprime divisor chain (see the discussion after Corollary \ref*{co01}).  Furthermore, there is no conclusion about the divisibility of LCM matrices by totally nonnegative GCD matrices yet. This motivates us to investigate on them.   Fortunately, the results of Guillot and Wu  in \cite{GWTN}  allow us to have the conclusion.  Moreover, we also introduce a \lq\lq column monotone matrix" from a given ordered set of monotone sequences in order to construct a totally nonnegative GCD matrix.

%--------------------------------------------\section{TN}----------------------------%--------------------------------------------\section{TN}----------------------------%--------------------------------------------\section{TN}----------------------------
\section{Totally nonnegative GCD matrices}
To be a self contained material, in this section, we collect the results of Guillot and Wu in \cite{GWTN} that will be used to prove our results. The necessary and sufficient conditions of totally nonnegative GCD matrices are given by:
\begin{thm}\label{TNt01} \cite{GWTN} Let $n\geq 3$ and let $S=\{x_1,\ldots,x_n\}$ be a set of distinct positive integers. Then the following statements are equivalent for the GCD matrix $(S)$: \begin{itemize}
		\item [(1)] $(S)$ is a totally nonnegative matrix.
		\item[(2)] For  $1\leq i\leq j\leq k\leq n$, $(x_{i},x_{k})=(x_{i},x_{j},x_{k})$ and $x_{j}\cdot(x_{i},x_{k})\mid x_{i}x_{k}$.
		\item[(3)] For $1\leq i\leq j\leq k\leq n$, $(x_i,x_j)\cdot(x_j,x_k)=x_{j}\cdot(x_i,x_k)$.
	\end{itemize}  Moreover, the above conditions imply that for all $1\leq i\leq j\leq k\leq l\leq n$,
	\begin{eqnarray}\label{TNe01}
	(x_i,x_k)\cdot(x_j,x_l)=(x_i,x_l)\cdot(x_j,x_k).
	\end{eqnarray}
\end{thm}
In particular case, we can rewrite the equation (\ref{TNe01}) as the following equation: 
\begin{eqnarray}\label{TNe02}
(x_{i},x_{j})=\frac{(x_{1},x_{j})\cdot(x_{i},x_{n})}{(x_{1},x_{n})},
\end{eqnarray}
for any  $1\leq i\leq j\leq n$.

Let $S=\{x_1,\ldots,x_n\}$ be a set of distinct positive integers and  let $2=p_1<p_2<p_3<\cdots$ denote the list of all prime numbers in increasing order.  For a given positive integer $x$, denote $e_t(x)\in \mathbb{N}_0$ (set of nonnegative integers) the power of $p_t$ occurring in the prime factorization of $x$ in which, by convention, $e_t(x)=0$ if $p_t$ is not a divisor of $x$.  Namely, for any $x\in \mathbb{N}$,
$$ x=\prod_{i=1}^{\infty}p_i^{e_i(x)}, $$
where only finitely many terms in the product are not equal to $1$.
An alternative necessary and sufficient condition for the totally nonnegative GCD matrix is given by:
\begin{thm}\label{TNt002}\cite{GWTN} Let $n\geq 3$ and let $S=\{x_1,\ldots,x_n\}$ be a set of distinct positive integers. Then the following statements are equivalent for the GCD matrix $(S)$:
	\begin{itemize}
		\item[(1)] $(S)$ is totally nonnegative.
		\item[(2)] For each $t\in \mathbb{N}$, the sequence $(e_t(x_i))_{i=1}^n$ is monotonic.
	\end{itemize}
\end{thm}
Here, a monotonic sequence means a sequence that is either non-decreasing or non-increasing.

To rewrite the above result, for a given ordered set of distinct positive integer $S=\{x_1,\ldots,x_n\}$,  we let $p_1<p_2<\cdots<p_k$ the list of all distinct prime factors of  $x:=\prod_{i=1}^nx_i$. For a given integer $i,j\in\mathbb{N}$, we denote by $a_{ij}$ the power of $p_{j}$  occurring in the prime factorization of $x_i$. Hence, for every $x_{i}\in S$, $$
	x_{i}=\prod\limits_{\substack{j=1}}^{k}p_{j}^{a_{ij}}.
$$
Denote by $$\operatorname{Pow}(S):=(a_{ij})\in M_{n\times k}(\mathbb{N}_0)$$ the $n\times k$ matrix constructed by the powers of prime factorization of $\prod_{x_i\in S}x_i$ where the primes are ordered by the usual order relation $<$.  For example, if $p_1<p_2<p_3$ are prime numbers and $$S=\{p_1^5 p_{3}^3,\,p_1^4 p_2 p_{3}^3,\,p_1^3 p_2 p_{3}^{2},\,p_1p_2^3,\,p_2^4\}\quad \hbox{ and }\quad S'=\{p_2^4,\, p_1^5 p_{3}^3,\,p_1^3 p_2 p_{3}^{2},\,p_1^4 p_2 p_{3}^3,\,\,p_1p_2^3\} $$  then, by the definition, \begin{equation}\label{ex1}
\operatorname{Pow}(S)=\left(\begin{array}{ccc}
	5&0&3\\
	4&1&3\\
	3&1&2\\
	1&3&0\\
	0&4&0
	\end{array}\right) \quad \hbox{ and } \quad \operatorname{Pow}(S')=\left(\begin{array}{ccc}
	0&4&0\\
	5&0&3\\
	3&1&2\\
	4&1&3\\
	1&3&0
	\end{array}\right) .
\end{equation} 
Note that $S$ and $S'$ are not gcd-closed sets and $S=S'$ as sets but not ordered sets.  
\begin{df}
Let $M\in\mathbb{M}_{n\times k}(\mathbb{R})$ be a matrix of size $n\times k$ over $\mathbb{R}$.  We say that $M$ is a \textbf{column monotone matrix} if each column of $M$ is a monotone sequence.
\end{df}
In (\ref{ex1}), the matrix $\operatorname{Pow}(S)$ is a column monotone matrix where as  $\operatorname{Pow}(S')$ is not. 
The following well known matrices are examples of column monotone matrices:
\begin{itemize}
		\item[$\bullet$]  \textbf{Vandermonde matrix}\\
	For $n$ real numbers $0<x_{1}<x_{2}<\cdots<x_{n}$, the vandermonde matrix $V(x_{1},x_{2},\ldots,x_{n})$ is defined to be 
	\begin{eqnarray*}
		V(x_{1},x_{2},\ldots,x_{n})=\left(\begin{array}{ccccc}
			1&x_{1}&x_{1}^{2}&\cdots&x_{1}^{n-1}\\
			1&x_{2}&x_{2}^{2}&\cdots&x_{2}^{n-1}\\
			\dot{:}&\dot{:}& &\dot{:}&\dot{:}\\
			1&x_{n}&x_{n}^{2}&\cdots&x_{n}^{n-1}
		\end{array}\right).
	\end{eqnarray*}
\item[$\bullet$]  \textbf{Hilbert matrix}\\
The Hilbert matrix $H_n=(h_{ij})$ is the  $n\times n$ matrix with 
 $h_{ij}=\frac{1}{i+j-1}$; namely
\begin{eqnarray*}	H_n=\left(\begin{array}{cccc} 
		1&\frac{1}{2}&\cdots&\frac{1}{n}\\
		\frac{1}{2}&\frac{1}{3}&\cdots&\frac{1}{n+1}\\
		\vdots&\vdots& \ddots&\vdots\\
		\frac{1}{n}&\frac{1}{n+1}&\cdots&\frac{1}{2n-1}
	\end{array}\right).
\end{eqnarray*}
	\item[$\bullet$] \textbf{Pascal matrix}  \\ The $n\times n$ matrix obtained by writing the Pascal's triangle as a symmetric matrix $P_{n}$ or a lower triangular matrix $L_{n}$.  For example, when $n=4$,
	\begin{eqnarray*}
		P_{4}=\left(\begin{array}{cccc}
			1&1&1&1\\
			1&2&3&4\\
			1&3&6&10\\
			1&4&10&20
		\end{array}\right)\text{ and }
	L_{4}=\left(\begin{array}{cccc}
		1&0&0&0\\
		1&1&0&0\\
		1&2&1&0\\
		1&3&3&1
	\end{array}\right).
\end{eqnarray*}
\end{itemize}
Note that for any prime numbers $p_{1}<p_{2}<p_{3}<p_4$,   the ordered set $$S=\{p_{1} p_{2} p_{3} p_{4},\,p_{1} p_{2}^2 p_{3}^3 p_{4}^4,\,p_{1}p_{2}^3 p_{3}^6 p_{4}^{10},\,p_{1} p_{2}^4 p_{3}^{10} p_{4}^{20}\}$$ has $Pow(S)=P_4$. By direct computation, we obtain that
	\begin{eqnarray}\label{pascalmatrix}
	(S)=\left(\begin{array}{cccc}
		p_{1} p_{2} p_{3} p_{4}&p_{1} p_{2} p_{3} p_{4}&p_{1} p_{2} p_{3} p_{4}&p_{1} p_{2} p_{3} p_{4}\\
		p_{1} p_{2} p_{3} p_{4}&p_{1} p_{2}^2 p_{3}^3 p_{4}^4&p_{1} p_{2}^2 p_{3}^3 p_{4}^4&p_{1} p_{2}^2 p_{3}^3 p_{4}^4\\
		p_{1} p_{2} p_{3} p_{4}&p_{1} p_{2}^2 p_{3}^3 p_{4}^4&p_{1} p_{2}^3 p_{3}^6 p_{4}^{10}&p_{1} p_{2}^3 p_{3}^6 p_{4}^{10}\\
		p_{1} p_{2} p_{3} p_{4}&p_{1} p_{2}^2 p_{3}^3 p_{4}^4&p_{1} p_{2}^3 p_{3}^6 p_{4}^{10}&p_{1} p_{2}^4 p_{3}^{10} p_{4}^{20}
	\end{array}\right).
\end{eqnarray}
The above GCD matrix $(S)$ in (\ref{pascalmatrix})  is constructed from the set $S$ that has $$\max_{x\in S}\{|{G_S(x)}|\}=1.$$
Then, by the result of Hong \cite{16}, $(S)\mid [S]$.  This conclusion can also be obtained by our main result.

Now by using our notations, we can rewrite Theorem \ref{TNt002} as: 
\begin{thm}(cf.\label{TNt02}\cite{GWTN}) Let $n\geq 3$ and let $S=\{x_1,\ldots,x_n\}$ be an ordered set of distinct positive integers. Then the following statements are equivalent for the GCD matrix $(S)$:
	\begin{itemize}
		\item[(1)] $(S)$ is totally nonnegative.
		\item[(2)] $\operatorname{Pow}(S)$ is a column monotone matrix. 
	\end{itemize}
\end{thm}

By Theorem \ref{TNt02}, in order to construct a totally nonnegative  GCD matrix $(S)$ of size $n\times n$, it can be done simply by picking a column monotone matrix $M:=(m_{ij})\in M_{n\times k}(\mathbb{N}_0)$ of size $n\times k$ over $\mathbb{N}_0$, for any $k\in \mathbb{N}$, and choosing primes $p_1<p_2<\cdots <p_k$.  Then, the ordered set $S:=\{x_1,x_2,\dots,x_n\}$, where $$x_i:=\prod_{j=1}^k p_j^{m_{ij}},$$ will serve the totally nonnegative GCD matrix $(S)$.

Another main tool for our calculations is the formula for the inverse of totally nonnegative GCD matrices.

\begin{thm}\label{TNt03}\cite{GWTN}
	Let $n\geq 3$ and let $S=\{x_1,\ldots,x_n\}$ be a set of distinct positive integers. Then the following statements are equivalent for the GCD matrix $(S)$:
	\begin{itemize}
		\item[(1)] $(S)$ is non-singular totally nonnegative matrix.
		\item[(2)] $(S)^{-1}$ is tridiagonal with nonzero superdiagonal elements.
		\end{itemize}
	\begin{eqnarray*}
		(S)^{-1}&=\left(\begin{array}{ccccc}
			b_1&a_2& & & \\
			a_2&b_2&a_3& & \\
			\,&\ddots&\ddots&\ddots& \\
			& & a_{n-1}&b_{n-1}&a_n\\
			& & &a_n&b_n
		\end{array}\right),
	\end{eqnarray*}
	where
	\begin{eqnarray*}
		a_{i+1}=\frac{(x_1,x_n)}{(x_i,x_n)(x_1,x_{i+1})-(x_{i+1},x_n)(x_1,x_i)}\,\,\,\,\,\,\,\,\text{if }1\leq i\leq n-1,
	\end{eqnarray*}
	and
	\begin{eqnarray*}
		b_{i}=\begin{cases}
			-\frac{(x_2,x_n)}{(x_1,x_n)}a_2&\text{if }i=1,\\
				-\frac{(x_{i-1},x_n)(x_1,x_{i+1})-(x_{i+1},x_n)(x_1,x_{i-1})}{(x_1,x_n)}a_ia_{i+1}&\text{if }2\leq i\leq n-1,\\
					-\frac{(x_1,x_{n-1})}{(x_1,x_n)}a_n&\text{if }i=n.
		\end{cases}
\end{eqnarray*}
\end{thm}

\section{Divisibility of [S] by (S)}
 As a result of  Beslin and  Ligh \cite{BL1989}, for any set of distinct positive integers $S$, $(S)$ is a positive definite matrix. This implies that $(S)^{-1}$ exists.   The following lemma is a main tool used to calculate the entries of $[S](S)^{-1}$.
\begin{lem}\label{lem301}
	Let $S = \{x_{1}, x_{2}, ..., x_{n}\}$ be a set of distinct positive integers such that $(S)$ is totally nonnegative. Then for any $1\leq i\leq j\leq n$,\begin{eqnarray}
	[x_{i},x_{j}]=\frac{(x_{1},x_{i})(x_{j},x_{n})}{(x_{1},x_{n})}.
	\end{eqnarray}
\end{lem}
\begin{proof}
	Let $S = \{x_{1}, x_{2}, \dots, x_{n}\}$ be a set of distinct positive integers such that $(S)$ is totally nonnegative. For any $1\leq i\leq j\leq n$, by (\ref{TNe02}), we have
	\begin{eqnarray*}
	[x_{i},x_{j}]&=& \frac{x_{i}\cdot x_{j}}{(x_{i},x_{j})}\\
	&=&\frac{(x_{i},x_{i})(x_{j},x_{j})}{(x_{i},x_{j})}\\
	&=&\frac{(x_{1},x_{i})(x_{i},x_{n})(x_{1},x_{j})(x_{j},x_{n})}{(x_{1},x_{n})(x_{1},x_{n})(x_{i},x_{j})}\\
	&=&\frac{(x_{1},x_{i})(x_{j},x_{n})}{(x_{1},x_{n})}\cdot\frac{(x_{1},x_{j})(x_{i},x_{n})}{(x_{1},x_{n})(x_{i},x_{j})}\\
	&=& \frac{(x_{1},x_{i})(x_{j},x_{n})}{(x_{1},x_{n})}\cdot\frac{(x_i,x_j)}{(x_i,x_j)}\\& =& \frac{(x_{1},x_{i})(x_{j},x_{n})}{(x_{1},x_{n})}.\\
	\end{eqnarray*}
\end{proof}
\begin{thm}\label{Dtn}
	Let $S = \{x_{1}, x_{2},\dots, x_{n}\}$ be a set of distinct positive integers. If $(S)$ is totally nonnegative, then \begin{eqnarray*}
		([S](S)^{-1})_{ij}=\begin{cases}
			-1&;\,i=j\neq1, n,\\
			\frac{x_{2}}{(x_{1},x_{2})}&;\,i=2,\,j=1,\\
			\frac{(x_{i},x_{n})}{(x_{1},x_{n})}&;\,i\neq1,2,\,j=1,\\
			\frac{x_{n-1}}{(x_{n-1},x_{n})}&;\,i=n-1,\,j=n,\\
			\frac{(x_{1},x_{i})}{(x_{1},x_{n})}&;\,i\neq n,n-1,\,j=n,\\
			0&; \text{ otherwise,}
		\end{cases}
	\end{eqnarray*}
	and $[S](S)^{-1}\in M_n(\mathbb{Z})$; namely, $(S)\mid[S]$.
\end{thm}
\begin{proof}
To simplify the notation, denote $U := [S](S)^{-1}.$ To verify the results, we will show that every entry $u_{ij}$ belong to  $\mathbb{N}$ by using Theorem \ref{TNt01}, Theorem \ref{TNt03} and Lemma \ref{lem301}.  We separate the calculation of  the entries of $U$ into 3 cases:  Case $j=1$, Case $j=n$ and Case $j\neq1,n$.

	In the first case, for $j=1$, we see from Theorem \ref{TNt03} that, for any $1\leq i \leq n$, 
	\begin{eqnarray*}
		U_{i1} &=& [x_{i}, x_{1}]b_{1} + [x_{i}, x_{2}]a_{2}\\
		&=&\Big(\frac{-[x_{i}, x_{1}](x_{2}, x_{n})}{(x_{1}, x_{n})} + [x_{i}, x_{2}] \Big)a_{2},		   
	\end{eqnarray*}
where $$b_1=-\frac{(x_2,x_n)}{(x_1,x_n)}a_2\; \hbox{ and }\; a_2=\frac{(x_1,x_n)}{(x_1,x_n)(x_1,x_{2})-(x_{2},x_n)(x_1,x_1)}.$$
	If $i=1$, then
	\begin{eqnarray*}
		U_{11} &=& \Big(\frac{-x_{1}(x_{2}, x_{n})}{(x_{1}, x_{n})} + [x_{1}, x_{2}] \Big)a_{2}\\
		&=& \Big(\frac{-(x_{1}, x_{1})(x_{2}, x_{n})}{(x_{1}, x_{n})} + \frac{(x_{1}, x_{1})(x_{2}, x_{n})}{(x_{1}, x_{n})} \Big)a_{2}\,\,\,\,\,\,\,\,\text{(by Lemma \ref{lem301})}\\
		&=&0.
	\end{eqnarray*}
	If $i=2$, then 
	\begin{eqnarray*}
		U_{21} &=& \Big(\frac{-[x_{2}, x_{1}](x_{2}, x_{n})}{(x_{1}, x_{n})} + x_{2} \Big)a_{2}\\
		&=& \Big(\frac{- x_{2}\cdot x_{1}(x_{2}, x_{n})}{(x_{2}, x_{1})(x_{1}, x_{n})} + x_{2} \Big)a_{2}\\
		&=& x_{2}\Big(\frac{-x_{1}(x_{2}, x_{n})+(x_{2}, x_{1})(x_{1}, x_{n})}{(x_{2}, x_{1})(x_{1}, x_{n})}\Big) \cdot \Big(\frac{(x_{1}, x_{n})}{(x_{1}, x_{n})(x_{1}, x_{2})-x_1(x_{2}, x_{n})} \Big)\\
		&=& \frac{x_{2}}{(x_{1}, x_{2})} \in \mathbb{N}.
	\end{eqnarray*}
If $i>2$, we have that
\begin{eqnarray*}
	U_{i1} &=& \Big(-\frac{[x_{i}, x_{1}](x_{2}, x_{n})}{(x_{1}, x_{n})} + [x_{i},x_{2}] \Big)a_{2}\\
	&=& \Big(-\frac{x_{1}\cdot(x_{i},x_{n})(x_{2},x_{n})}{(x_{1}, x_{n})(x_{1}, x_{n})} + \frac{(x_{i},x_{n})(x_{1},x_{2})}{(x_{1},x_{n})} \Big)a_{2}\,\,\,\,\,\,\,\,\text{(by Lemma \ref{lem301})}\\
	&=& (x_{i},x_{n})\Big(\frac{-x_{1}\cdot(x_{2},x_{n})+(x_{1},x_{2})(x_{1},x_{n})}{(x_{1}, x_{n})(x_{1}, x_{n})}\Big)\cdot\Big( \frac{(x_{1},x_{n})}{(x_{1},x_{2})(x_{1},x_{n})-x_1\cdot(x_{2},x_{n})} \Big)\\
	&=& \frac{(x_{i}.x_{n})}{(x_{1}, x_{n})}.
\end{eqnarray*}
By Theorem \ref{TNt01}(2), $(x_1,x_n)=(x_1,x_2,\dots,x_n)$.  Then $(x_1,x_n)\mid (x_i,x_n)$; i.e., $ U_{i1} \in \mathbb{N}$ for all $i>2$.

In the second case, for $j=n$, we see from Theorem \ref{TNt03} that, for any $1\leq i \leq n$, 
\begin{eqnarray*}
	U_{in} &=& [x_{i}, x_{n-1}]a_{n} + [x_{i}, x_{n}]b_{n}\\
	&=&[x_{i}, x_{n-1}]a_{n} - \frac{[x_{i}, x_{n}](x_{1},x_{n-1})}{(x_{1},x_{n})}a_{n}\\ 
	&=& \Big([x_{i}, x_{n-1}] - \frac{[x_{i}, x_{n}](x_{1},x_{n-1})}{(x_{1},x_{n})}\Big)a_{n},
\end{eqnarray*}
where $$	a_{n}=\frac{(x_1,x_n)}{(x_{n-1},x_n)(x_1,x_{n})-(x_{n},x_n)(x_1,x_{n-1})} \;\;\hbox{ and } \;\; b_n=-\frac{(x_1,x_{n-1})}{(x_1,x_n)}a_n.$$
If $i=n$, we have that
\begin{eqnarray*}
U_{nn}&=&\Big([x_{n},x_{n-1}]-\frac{[x_{n},x_{n}](x_{1},x_{n-1})}{(x_{1},x_{n})}\Big)a_{n}\\
&=&\Big(\frac{(x_{1},x_{n-1})(x_n,x_n)}{(x_{1},x_{n})}-\frac{[x_{n},x_{n}](x_{1},x_{n-1})}{(x_{1},x_{n})}\Big)a_{n}\,\,\,\,\,\,\,\text{(by Lemma \ref{lem301})}\\
&=&0.
\end{eqnarray*}
If $i=n-1$, we have that
\begin{eqnarray*}
	U_{n-1,n}&=&\Big([x_{n-1},x_{n-1}]-\frac{[x_{n-1},x_{n}](x_{1},x_{n-1})}{(x_{1},x_{n})}\Big)a_{n}\\
	&=&\Big(x_{n-1}-\frac{x_{n-1}\cdot x_{n}\cdot(x_{1},x_{n-1})}{(x_{n-1},x_{n})(x_{1},x_{n})}\Big)a_{n}\\
		&=&x_{n-1}\cdot\Big(\frac{(x_{n-1},x_{n})(x_{1},x_{n})-(x_{n},x_{n})(x_{1},x_{n-1})}{(x_{n-1},x_{n})(x_{1},x_{n})}\Big)\cdot\Big(\frac{(x_{1},x_{n})}{(x_{n-1},x_{n})(x_1,x_{n})-(x_{n},x_{n})(x_{1},x_{n-1})}\Big)\\
		&=& \frac{x_{n-1}}{(x_{n-1},x_{n})} \in \mathbb{N}.
\end{eqnarray*}
If $i<n-1$, we have
\begin{eqnarray*}
	U_{in} &=& \Big([x_{i}, x_{n-1}] - \frac{[x_{i}, x_{n}](x_{1},x_{n-1})}{(x_{1},x_{n})}\Big)a_{n}\\
	&=& \Big(\frac{(x_{1},x_{i})(x_{n-1},x_{n})}{(x_{1},x_{n})}-\frac{(x_{1},x_{i})(x_{n},x_{n})(x_{1},x_{n-1})}{(x_{1},x_{n})(x_{1},x_{n})}\Big)a_{n}\,\,\,\,\,\,\text{(by Lemma \ref{lem301})}\\
	&=& (x_{1},x_{i})\cdot\Big(\frac{(x_{n-1},x_{n})(x_{1},x_{n})-(x_{n},x_{n})(x_{1},x_{n-1})}{(x_{1},x_{n})(x_{1},x_{n})}\Big)\cdot\Big(\frac{(x_{1},x_{n})}{(x_{n-1},x_{n})(x_1,x_{n})-(x_{n},x_{n})(x_{1},x_{n-1})}\Big)\\
	&=& \frac{(x_{1},x_{i})}{(x_{1},x_{n})}.
\end{eqnarray*}
By Theorem \ref{TNt01}(2), $(x_1,x_n)=(x_1,x_2,\dots,x_n)$.  Then $(x_1,x_n)\mid (x_1,x_i)$; i.e., $ U_{in} \in \mathbb{N}$ for all $1\leq i<n-1$.

Finally, for the case $j\neq1,n$, by using Theorem \ref{TNt03}, we have that, for $i=1,\dots,n$
\begin{eqnarray}\label{finalcase}
U_{ij}=[x_{i},x_{j-1}]a_{j}+[x_{i},x_{j}]b_{j}+[x_{i},x_{j+1}]a_{j+1},
\end{eqnarray}
where $$
	a_{j+1}=\frac{(x_1,x_n)}{(x_j,x_n)(x_1,x_{j+1})-(x_{j+1},x_n)(x_1,x_j)}$$ and $$
b_j=-\frac{(x_{j-1},x_n)(x_1,x_{j+1})-(x_{j+1},x_n)(x_1,x_{j-1})}{(x_1,x_n)}a_ja_{j+1}.
$$
If $i=j$, by substituting the value of $a_i,b_i,a_{i+1}$ into the equation (\ref{finalcase}) and aligning the equation, then 
\begin{eqnarray*}
U_{ii}=\frac{A+B+C}{\Big((x_{i-1},x_{n})(x_1,x_{i})-(x_{i},x_{n})(x_1,x_{i-1})\Big)\Big((x_{i},x_{n})(x_{1},x_{i+1})-(x_{i+1},x_{n})(x_{1},x_{i})\Big)}(x_1,x_{n}),
\end{eqnarray*}
where
\begin{eqnarray*}
A&=&[x_{i},x_{i-1}]\Big((x_{i},x_{n})(x_{1},x_{i+1})-(x_{i+1},x_{n})(x_{1},x_{i})\Big),\\
B&=&-[x_{i},x_{i}]\Big((x_{i-1},x_{n})(x_{1},x_{i+1})-(x_{i+1},x_{n})(x_{1},x_{i-1})\Big),\\
C&=&[x_{i},x_{i+1}]\Big((x_{i-1},x_{n})(x_1,x_{i})-(x_{i},x_{n})(x_1,x_{i-1})\Big).
\end{eqnarray*}
We now rewrite $A$ as: 
\begin{eqnarray*}
A&=&[x_{i},x_{i-1}]\Big((x_{i},x_{n})(x_{1},x_{i+1})-(x_{i+1},x_{n})(x_{1},x_{i})\Big)\\
&=& \frac{(x_{1},x_{i})(x_{i},x_{n})(x_{1},x_{i-1})(x_{i-1},x_n)}{(x_i,x_{i-1})(x_{1},x_{n})(x_{1},x_n)}\cdot\Big((x_{i},x_{n})(x_{1},x_{i+1})-(x_{i+1},x_{n})(x_{1},x_{i})\Big)\\
&=&\frac{(x_{1},x_{i})(x_{i-1},x_n)(x_{i},x_{n})(x_{1},x_{i+1})}{(x_i,x_{i-1})(x_{1},x_{n})(x_{1},x_n)}-\frac{(x_{1},x_{i})(x_{i-1},x_n)(x_{i+1},x_{n})(x_{1},x_{i})}{(x_i,x_{i-1})(x_{1},x_{n})(x_{1},x_n)}\,\,\, \text{(by Lemma \ref{lem301}})\\
&=&\frac{(x_i,x_n)(x_1,x_{i-1})(x_i,x_n)(x_1,x_{i+1})-(x_i,x_{n})(x_1,x_{i-1})(x_{i+1},x_n)(x_{i},x_{1})}{(x_{1},x_{n})}\,\,\, \text{(by Lemma \ref{lem301}}).
\end{eqnarray*}
Similarly, we can calculate that
\begin{eqnarray*}
B&=& \frac{(x_{i+1},x_n)(x_1,x_{i-1})(x_1,x_i)(x_i,x_n)-(x_{i-1},x_{n})(x_1,x_{i+1})(x_{1},x_i)(x_{i},x_{n})}{(x_{1},x_{n})}\\
C&=&\frac{(x_1,x_i)(x_{i+1},x_{n})(x_{i-1},x_n)(x_1,x_{i})-(x_1,x_{i})(x_{i+1},x_{n})(x_{1},x_n)(x_{1},x_{i-1})}{(x_{1},x_{n})}.
\end{eqnarray*}
Obviously, 
\begin{eqnarray*}
A+B+C=-\frac{\Big((x_{i-1},x_{n})(x_1,x_{i})-(x_{i},x_{n})(x_1,x_{i-1})\Big)\Big((x_{i},x_{n})(x_{1},x_{i+1})-(x_{i+1},x_{n})(x_{1},x_{i})\Big)}{(x_{1},x_{n})}
\end{eqnarray*}
which implies that $U_{ii}=-1$.

If $i\neq j$, then we will consider separately two cases: Case $i> j$ and Case $i<j$.  For $i>j$, we note that $i\geq j+1>j-1$. Then, by substituting the value of $a_j,b_j,a_{j+1}$ into the equation (\ref{finalcase}) and aligning the equation, we have that
\begin{eqnarray*}
	U_{ij}
	&=&\frac{X+Y+Z}{\Big((x_{j-1},x_{n})(x_1,x_{j})-(x_{j},x_{n})(x_1,x_{j-1})\Big)\Big((x_{j},x_{n})(x_{1},x_{j+1})-(x_{j+1},x_{n})(x_{1},x_{j})\Big)}(x_1,x_{n}),
\end{eqnarray*}
where
\begin{eqnarray*}
X&=&[x_{i},x_{j-1}]\Big((x_{j},x_{n})(x_{1},x_{j+1})-(x_{j+1},x_{n})(x_{1},x_{j})\Big)\\
Y&=&-[x_{i},x_{j}]\Big((x_{j-1},x_{n})(x_{1},x_{j+1})-(x_{j+1},x_{n})(x_{1},x_{j-1})\Big)\\
Z&=&[x_{i},x_{j+1}]\Big((x_{j-1},x_{n})(x_1,x_{j})-(x_{j},x_{n})(x_1,x_{j-1})\Big).
\end{eqnarray*} In the same spirit as we have computed $A, B$ and $C$, we calculate that
\begin{eqnarray*}
X
&=&\frac{(x_i,x_n)(x_1,x_{j-1})(x_{j},x_{n})(x_{1},x_{j+1})-(x_i,x_n)(x_1,x_j)(x_{j+1},x_n)(x_{1},x_{j-1})}{(x_1,x_n)},\\
Y
&=& \frac{(x_i,x_{n})(x_{1},x_{j})(x_{j+1},x_{n})(x_{1},x_{j-1})-(x_i,x_n)(x_1,x_{j})(x_{j-1},x_{n})(x_{1},x_{j+1})}{(x_1,x_n)},\\
Z&=&\frac{(x_i,x_n)(x_1,x_{j})(x_{j-1},x_{n})(x_{1},x_{j+1})-(x_i,x_n)(x_1,x_{j-1})(x_{j},x_{n})(x_{1},x_{j+1})}{(x_1,x_n)}.
\end{eqnarray*}
It is clear that $X+Y+Z=0$, which implies that $U_{ij}=0$ for all $i>j$.  Next, if $i<j$, then $j+1>j-1\geq i$. In the same way as in the case $i>j$, we obtain that $U_{ij}=0$. Then the matrix $U$ can be expressed as: \begin{equation*}
 U=\left(\begin{array}{ccccccccc}
0&0&0&0&0&\cdots&0&0&\frac{x_{1}}{(x_{1},x_{n})}\\
\frac{x_{2}}{(x_{1},x_{2})}&-1&0&0&0&\cdots&0&0&\frac{(x_{1},x_{2})}{(x_{1},x_{n})}\\
\frac{(x_{3},x_{n})}{(x_{1},x_{n})}&0&-1&0&0&\cdots &0&0&\frac{(x_{1},x_{3})}{(x_{1},x_{n})}\\
\frac{(x_{4},x_{n})}{(x_{1},x_{n})}&0&0&-1&0&\cdots &0&0&\frac{(x_{1},x_{4})}{(x_{1},x_{n})}\\
\frac{(x_{5},x_{n})}{(x_{1},x_{n})}&0&0&0&-1&\cdots &0&0&\frac{(x_{1},x_{5})}{(x_{1},x_{n})}\\
\vdots&\vdots&\vdots&\vdots&\vdots&\ddots&\vdots&\vdots&\vdots\\
\frac{(x_{n-1},x_{n})}{(x_{1},x_{n})}&0&0&0&0&\cdots &-1&0&\frac{x_{n-1}}{(x_{1},x_{n})}\\
\frac{x_{n}}{(x_{1},x_{n})}&0&0&0&0&\cdots &0&0&0\\
\end{array}\right).
\end{equation*}
Since each entry of $U$ belongs to $\mathbb{Z}$, we conclude that $(S)\mid[S]$.
\end{proof}

It is well known that the divisibility of LCM matrices by GCD matrices constructed from $S$ is independent on the rearrangement of the entries of $S$.  Namely, for any set $S=\{x_1,x_2,\dots,x_n\}\subseteq \mathbb{N}$ with $|S|=n$, if $(S)\mid [S]$, then $(\sigma(S))$, where $$\sigma(S):=\{x_{\sigma(1)},x_{\sigma_2},\dots,x_{\sigma(n)}\},$$
still divides $[\sigma(S)]$ for any $\sigma\in S_n$ (set of all permutations of degree $n$).  Then, by Theorem \ref{Dtn} and Theorem \ref{TNt02}, we have that:  \begin{thm}\label{Dtn2}
	Let $S = \{x_{1}, x_{2}, \dots, x_{n}\}$ be an ordered set of distinct positive integers. If there exist $\sigma\in S_n$ such that $\operatorname{Pow}(\sigma(S))$ is a column monotone matrix, then $(S)\mid [S]$.
\end{thm} 
In (\ref{ex1}), the elements in ordered set $S'$ can be rearranged to the ordered set $S$ in which $\operatorname{Pow}(S)$ is a column monotone matrix. Then, by Theorem \ref{Dtn2}, $(S)\mid[S]$ and thus $(S')|[S']$.  The GCD matrix $(S)$ in (\ref{pascalmatrix})  is also a divisor of $[S]$ since $\operatorname{Pow}(S)=P_4$ is a column monotone matrix.

For a given totally ordered set $S$, a kernel $K:S\times S\rightarrow\mathbb{R}$ is said to be \textit{totally nonnegative} if $(K(x_{i},x_{j}))_{i,j=1}^{n}$ is totally nonnegative matrix for any choice of integers $x_{1}<\cdots<x_{n}$ in $S$ and any $n\geq1$ \cite{K22}. Then, the greatest common divisor function, $\operatorname{gcd}:S\times S\longrightarrow \mathbb{N}$, is a totally nonnegative kernel on $S$, when $\operatorname{Pow}(S)$ is a column monotone matrix.  Guillot and Wu \cite{GWTN}, proposed a property of a function $f:\mathbb{N}\longrightarrow \mathbb{R}$ preserving totally nonegative kernel $K$ on $S\subseteq \mathbb{N}$; namely, if $f$ is a multiplicative function and $f(x)\leq f(y)$ for every $x,y\in\mathbb{N}$ such that $x\mid y$, then $f\circ K$ is a totally nonnegative on $S$.  In particular, for any $e\in \mathbb{N}$, the e-th power function $\xi_e:\mathbb{N}\longrightarrow \mathbb{R}$, $\xi_e(x)=x^e$ for all $x\in S$, satisfies the conditions and thus $\xi_e\circ K$ is totally nonnegative for any $e\in \mathbb{N}$ and totally nonnegative kernel $K:S\times S\longrightarrow \mathbb{N}$.  It is also well known that $(\xi_e(S)):=(S^e)$ is a GCD matrix and $[\xi_e(S)]:=[S^e]$ is an LCM matrix. So, by Theorem \ref{TNt02} and Theorem \ref{Dtn2}, we have that:

\begin{cor}\label{co01}
Let $S = \{x_{1}, x_{2}, \dots, x_{n}\}$ be an ordered set of distinct positive integers. If there exists $\sigma\in S_n$ such that $\operatorname{Pow}(\sigma(S))$ is a column monotone matrix, then, for any $e\in \mathbb{N}$, $(S^e)$ is a totally nonnegative matrix and $(S^e)\mid [S^e]$.
\end{cor}

 A set of distinct positive integers $S = \{x_{1}, x_{2}, \dots, x_{n}\}$ is called a \textit{divisor chain} if there exists a permutation $\sigma\in S_n$ such that $x_{\sigma(1)}\mid x_{\sigma(2)}\mid\cdots\mid x_{\sigma(n)}$. The set $S$ is called a \textit{finitely many coprime divisor chains} if there is a positive integer $k$ such that $S$ can be partitioned as $S=S_{1}\cup S_{2}\cup\cdots\cup S_{k}$, where $S_{i}$ are divisor chains for all $1\leq i\leq k$ and each element in $S_{i}$ is coprime to each element in $S_{j}$ for all $1\leq i<j\leq k$. In 2020, Tan and Lin \cite{33} provided that for any finitely many coprime set $S$, $(S^e)\mid[S^e]$ for any positive integers $e$. However, for prime numbers $p_1<p_2<p_3< p_4$, the set  $$S=\{ p_2^9p_4^5,\; p_2^8 p_3 p_4^5,\;p_1 p_2^7 p_3 p_4^3,\;p_1 p_2^5 p_3^3 p_4^2,\;p_1^5 p_2^2 p_3^8 p_4^2,\;p_1^7 p_2 p_3^{11}\}$$
 is not a gcd-closed set nor a finitely many coprime divisor chains, since\begin{center}
	$(
p_1 p_2^5 p_3^3 p_4^2,p_1^5 p_2^2 p_3^8 p_4^2)=
p_1 p_2^{2} p_3^{3} p_4^2\notin S$
\end{center} and $S$ is separated into $6$ distinct chains but these chains are not coprimes.  Fortunately, \begin{center}
$\operatorname{Pow}(S)=\left(\begin{array}{cccc}
0&9&0&5 \\
0&8&1& 5 \\
1&7&1&3 \\
1&5 &3&2\\
5&2&8&2\\
7&1 &11 &0
\end{array}\right)$
\end{center} is a column monotone matrix. By Corollary \ref{co01}, we can conclude that $(S^e)\mid[S^e]$ for any positive integer $e$.

\section*{Acknowledgments}
The second author would like to thank Faculty of Science, Naresuan University, Phitsanulok, Thailand, for financial support on the project number P2563C102.

\bigskip

\address \textbf{Peeraphat Gatephan} 

{ Department of Mathematics, Faculty of Science, \\ Naresuan University, Phitsanulok 65000, Thailand}\\
\email{peeraphatg@nu.ac.th, \quad \quad}
\bigskip

\address \textbf{Kijti Rodtes} 

{ Department of Mathematics, Faculty of Science, \\ Naresuan University, Phitsanulok 65000, Thailand}\\
\email{kijtir@nu.ac.th, \quad \quad}


\begin{thebibliography}{20}
	
\bibitem{00} E. Altini\c{s}ik, M. Yildiz, A. Keskin, Non-divisibility of LCM matrices by GCD matrices on gcd-closed sets, Linear Algebra Appl. 516 (2017) 47-68.

\bibitem{BL1989} S. Beslin, S. Ligh, Greatest common divisor matrices, Linear Algebra Appl. 118 (1989), 69-76.

%\bibitem{BLA} S. Beslin, S. Ligh, Another generalization of Smith's determinant, Bull. Austral. Math. Soc. (3) 40 (1989), 413-415.

\bibitem{04}K. Bourque, S. Ligh, On GCD and LCM matrices, Linear Algebra Appl. 174 (1992) 65-74.

%\bibitem{05} K. Bourque, S. Ligh, Matrices associated with arithmetical functions, Linear Multilinear Algebra 34 (1993) 261–267.

%\bibitem{06}K. Bourque, S. Ligh, Matrices associated with multiplicative functions, Linear Algebra Appl. 216 (1995) 267–275.

\bibitem{07} W. Feng, S. Hong, J. Zhao, Divisibility properties of power LCM matrices by power GCD matrices on gcd-closed sets, Discrete Math. 309 (2009) 2627-2639.

\bibitem{GWTN} D. Guillot, J. Wu, Total nonnegativity if GCD matrices and kernels, Linear Algebra Appl. 578 (2019) 446-461.

%\bibitem{08} P. Haukkanen, I. Korkee, Notes on the divisibility of GCD and LCM matrices, Int. J. Math. Math. Sci. 6 (2005) 925–935.

\bibitem{13}S. Hong, On the factorization of LCM matrices on gcd-closed sets, Linear Algebra Appl. 345 (2002) 225-233.

\bibitem{16}S. Hong, Nonsingularity of matrices associated with classes of arithmetical functions on lcm-closed sets, Linear Algebra Appl. 416 (2006) 124-134.

%\bibitem{17}S. Hong, Divisibility properties of power GCD matrices and power LCM matrices, Linear Algebra Appl. 428 (2008) 1001–1008.

%\bibitem{21}S. Hong, J. Zhao, Y. Yin, Divisibility properties of Smith matrices, Acta Arith. 132 (2008) 161–175.

%\bibitem{22}L.-K. Hua, Introduction to Number Theory, Springer, New York, 1982.
\bibitem{K22} S. Karlin, Total Positivity, vol. 1, Stanford University ress, 1968.
\bibitem{S1876}H.J.S. Smith, On the value of a certain arithmeticaldeterminant, Proc. Lond. Math. Soc. 7 (1876) 208-212.

\bibitem{30}Q. Tan, Notes on non-divisibility of determinants of power GCD matrices and power LCM matrices, Southeast Asian Bull. Math. 33 (2009) 563-567.

%\bibitem{31}Q. Tan, Divisibility among power GCD matrices and among power LCM matrices on two coprime divisor chains, Linear Multilinear Algebra 58 (2010) 659–671.

%\bibitem{32} Q. Tan, M. Li, Divisibility among power GCD matrices and among power LCM matrices on finitely many coprime divisor chains, Linear Algebra Appl. 438 (2013) 1454–1466.

\bibitem{33}Q. Tan, Z. Lin, Divisibility of determinants of power GCD matrices and power LCM matrices on finitely many quasi-coprime divisor chains, Appl. Math. Comput. 217 (2010) 3910-3915.

\bibitem{TLC2016}Q. Tan, Z. Lin, L. Chen, Divisibility among power matrices associated with arithematic functions on finitely many quasi-coprime divisor chains, Linear Multilinear Algebra 64 (2016) 2030-2048.

%\bibitem{34}Q. Tan, Z. Lin, L. Liu, Divisibility among power GCD matrices and among power LCM matrices on two coprime divisor chains II, Linear Multilinear Algebra 59 (2011) 969–983.

\bibitem{35}Q. Tan, M. Luo, Z. Lin, Determinants and divisibility of power GCD and power LCM matrices on finitely many coprime divisor chains, Appl. Math. Comput. 219 (2013) 8112-8120.

\bibitem{37}J. Zhao, A characterization of the gcd-closed set S with $|S| =4$ such that $(S^e)$ divides $[S^e]$, J. Sichuan Univ. Nat. Sci. Ed. 45 (2008) 475–477.

\bibitem{38}W. Zhao, J. Zhao, A characterization for the gcd-closed set S with $|S| =5$ such that $(S^e)$ divides $[S^e]$, Southeast Asian Bull. Math. 33 (2009) 1023-1028.

%\bibitem{39}W. Zhao, Divisibility of power LCM matrices by power GCD matrices on gcd-closed sets, Linear Multilinear Algebra 62 (2014) 735–748.	
		
\end{thebibliography}
\end{document}